\pdfoutput=1
\documentclass{amsart}

\usepackage{graphicx}
\usepackage{dsfont}

\newtheorem{thm}{Theorem}[section]
\newtheorem{cor}[thm]{Corollary}
\newtheorem{lem}[thm]{Lemma}

\theoremstyle{definition}

\theoremstyle{remark}
\newtheorem{rem}[thm]{Remark}
\numberwithin{equation}{section}

\begin{document}

\title[]{Palindromic Periodicities}

\author{Jamie Simpson}%
\address{130 Preston Point Rd \newline
\indent East Fremantle,\newline
\indent WA 6845\newline
\indent Australia \newline}%
\email{jamiesimpson320@gmail.com}

\begin{abstract} A palindromic periodicity is a factor of an infinite word $(ps)^\omega$ where $p$ and $s$ are palindromes and the factor has length at least $|ps|$, for example $accabaccab$.   In this paper we describe several ways in which a palindromic periodicity may arise through the interaction of palindromes and periodicity, the simplest case being when a palindrome is itself periodic. We then consider what happens when a word is a palindromic periodicity in two ways, a situation similar to that considered in the Fine and Wilf Lemma and obtain something slightly stronger than that lemma. The paper ends with suggestions for further work.\\

\noindent \emph{Keywords:} palindrome, periodicity, word.
\end{abstract}

\maketitle

\section{Introduction}
We are concerned with words which are factors of infinite words of the form $(ps)^\omega$ where $p$ and $s$ are palindromes and which have length at least $|ps|$. We call such words \emph{palindromic periodicities}. The paper is set out as follows. After reviewing notation we discuss some of the simple properties of these words.  We then describe, in Section 2, several ways in which palindromic periodicities can arise. In the third section we consider ways in which these words interact and obtain an analogy of Fine and Wilf's Periodicity Lemma for these words.  Fine and Wilf does, of course, already apply to our periodic words but their extra structure means a slightly stronger result is possible.   We also obtain a  result about a palindrome embedded in a palindromic periodicity.  We end the paper with some discussion and suggestions for further investigation.

 We use the usual notation for combinatorics of words. A word containing $n$
letters is $w=w[1 \dots n]$, with $w[i]$ being the $i$th letter and
$w[i\dots j]$ the \emph{factor} beginning at position $i$ and ending at
position $j$.  If $i=1$ then the factor is a \emph{prefix} and if
$j=n$ it is a \emph{suffix}.  If a prefix and a suffix of $w$ both equal a word $x$ then $x$ is a \emph{border} of $w$. If $w=tur$ where $t$, $u$ and $v$ are factors, we say that $w$ is the \emph{union} of $tu$ and $ur$. A factor which is neither a suffix nor
a prefix is a \emph{proper factor}\footnote{In linguistics a proper factor is called an ``infix". Linguists use $\emptyset$ to mean empty set though they don't call it an empty set  and iff as we do but they pronounce it, at least in the Linguistics Department of the University of Western Australia, as ``if if".}. The \emph{length} of $w$, written $|w|$, is the
number of letters that $x$ contains.  If $w=uv$ where $u$ and $v$ are words then $vu$ is a
\emph{conjugate} of $w$.  The \emph{empty word} $\epsilon$ is a word
with length 0.  A word or factor $w$ is \emph{periodic} with
period $p$ if $w[i]=w[i+p]$ for all $i$ such that both $w[i]$ and
$w[i+p]$ are in $w$. $w^{\omega}=www\dots$ is the infinite concatenation of $w$ with itself

\begin{rem} \label{border implies periodic}
It is easy to see that if $x$ is a border of $w$ then $w$  has period $|w|-|x|$.
\end{rem}
 We will use the following well-known
propositions.

\begin{lem} \label{FW} \textup{(Periodicity Lemma of Fine and Wilf, \cite{FW})} Let $w$
be a word having two periods $p$ and $q$.  If $|w| \ge
p+q-\gcd(p,q)$ then $w$ also has period $\gcd(p,q)$.
\end{lem}

\begin{lem} \label{LA} (Lemma 2.1 of \cite{CMR} and Lemma 8.1.1 of \cite{Lo2}) Let $w$ be a word having two periods $p$ and
$q$ with $p>q$. Then the suffix and prefix  of length $|w|-q$ both
have period $p-q$.\end{lem}

\begin{lem} \label{L4} \textup{(Lemma 8.1.3 of \cite{Lo2})} Let $w$
be a word with period $q$ which has a factor $u$ with $|u| \ge q$
that has period $r$, where $r$ divides $q$. Then $w$ has period $r$.
\end{lem}

\begin{lem} \label{overlapping periodicities}\textup{(Lemma 8.1.2 of \cite{Lo2})}
Let  $u$, $v$ and $w$ be words such that $uv$ and $vw$ both have
period $p$ and $|v| \ge p.$ Then the word $uvw$ has period $p$.
\end{lem}

\begin{lem} \label{new lemma}
  If a word $w$ has period $p$ and $w[i+1..i+q] =w[j+1..j+q]$ where $i+1<j<i+q$ and $q\ge p$ then $w$ has period $\gcd(p,j-i)$,
  \begin{proof}
    The factor $w[i+1..j+q]$ has a border of length $q$ and therefore, by Remark \ref{border implies periodic}, has period $j-i$. It also has period $p$ and length $j+q-i$.  Since $q \ge p$ the periodicity lemma applies and the factor has period $\gcd(j-i,p)$. By Lemma \ref{L4} this periodicity extends to the whole of $w$.
  \end{proof}
\end{lem}

The \emph{reverse} of a word $w[1\dots n]$ is the word
$R(w)=w[n]w[n-1]\dots w[1]$. The word $w$ is a \emph{palindrome} if
 $w=R(w)$. Thus the empty word $\epsilon$ is a palindrome.  A palindrome is \emph{odd} or
\emph{even} if its length is, respectively, odd or even.  If
$w[i\dots j]$ is a palindrome we say it has \emph{centre} $c=(i+j)/2$
and \emph{radius} $r=(j-i)/2$. Note that $c$ and $r$ are integers if
the palindrome is odd and each is an integer plus 1/2 if the
palindrome is even. Much of what follows concerns centres of palindromes and it will be useful to have the following notation:
$$\mathds{Z}/2=\{n/2:n \in \mathds{Z}\}$$ so that $c$ and $r$  are in $\mathds{Z}/2$ whereas $2c$, $2r$ and $c+r$ are in $\mathds{Z}$ .
If $w$ contains a palindrome with centre $c$ and radius $r$ and $c-r \le i<j \le c+r$ then
\begin{equation}\label{e1}
w[i]=w[2c-i]
\end{equation} and
\begin{equation}\label{e1.5}
w[i\dots j]=w[2c-j..2c-i].
\end{equation}
Having centres which are not integers means we have to decide whether $w[5/2]$ is in $w[1..2]$ or in $w[3..4]$, or in both or in neither. ``Neither" seems wrong since it is in the concatenation of the two factors and ``both" would also be inconvenient. We therefore adopt the convention that if $a$ and $b$ are integers then $w[a-1/2]$ is in $w[a..b]$ but $w[b+1/2]$ is not\footnote{We could adopt a more elaborate notation here, say with both $w[1/2]$ and $w[5/2]$ in $w(1..2)$ and neither in $w[1..2]$, but this is different from the notation for open and closed intervals where $(1,2)$ is contained in $[1,2]$.}.

 If $w[i\dots j]$ is a palindrome with
$j \ge i+2$ then so is $w[i+1\dots j-1]$. If $w[c-r\dots c+r]$ is a
palindrome but none of $w[c-r-1\dots c+r+1]$, $w[c-r\dots c+r+1]$, and
$w[c-r-1\dots c+r]$ is, then we say $w[c-r\dots c+r]$ is a \emph{maximal
palindrome}. The second and third cases here mean that we do not,
for example, consider $aa$ to be maximal in $baaac$. If
$w[c-r\dots c+r]$ is maximal then the palindromes $w[c-r+i\dots c+r-i]$,
$i=1,\dots,\lfloor r \rfloor$, are \emph{nested} in $w[c-r\dots c+r]$.
If a palindrome is even, respectively odd, then its nested
palindromes are even, respectively odd.


We now make some observations about palindromic periodicities.  Recall that
a palindromic periodicity is a factor of length at least $|ps|$ of the infinite word $(ps)^{\omega}$ where $p$ and $s$ are palindromes.  For example $baccaba$ is a factor of $(ps)^\omega$ where $p=aba$ and $s=cc$. It is also a factor of $(ps)^{\omega}$ when $p=b$ and $s=acca$, or when $p=cabac$ and $s=\epsilon$.  In fact $ps$ can always be replaced with a single palindrome except when both $p$ and $s$ are odd palindromes.  From similar considerations we  make the following observation.

\begin{rem} \label{remark} Any palindromic periodicity is a prefix of some infinite word $(ps)^{\omega}$ where $p$ and $s$ are palindromes.
\end{rem}
Thus, while $p$ and $s$ are not specified for a given palindromic periodicity, their centres are fixed, and the distance between these centres is fixed and equals $|ps|/2$, which we call the \emph{half-period} of the palindromic periodicity. This will sometimes be a more convenient parameter than the full period.  Note that the half-period is in $\mathds{Z}/2$. The palindromic periodicity $abbcbbadabbcbbdb$  is a factor of $(ps)^{\omega}$ with $p=ada$ and $s=bbcbb$ with centres $d$ and $c$. We call such centres \emph{essential palindromic centres}, or just \emph{essential centres}.  This palindromic periodicity also contains the palindrome $bb$ but its centre is not the centre of $p$ or $s$ and it is not essential.  We note the following:
\begin{rem} Any essential centre of a palindromic periodicity $w$ is the centre of a palindromic prefix or a palindromic suffix of $w$. This is not  true of non-essential centres.
\end{rem}

Following these observations we see that the following is an alternative definition of a palindromic periodicity:  word  $w$ is a \emph{palindromic periodicity} with length $n$, offset $r$ and half-period $h$ if each position $r, r+h,\dots r+\lfloor{(n-r)}/h\rfloor h $ is the palindromic centre of a palindromic prefix or suffix, and further, that the sum of the lengths of the longest such prefix and the longest such suffix is at least $|w|$. The centres of the longest palindromic prefix and longest palindromic suffix are, respectively

$$r+ \left \lfloor \frac{\lfloor \frac{n-r}{h}\rfloor}{2}\right \rfloor h \text{ and } r+ \left \lceil \frac{\lfloor \frac{n-r}{h}\rfloor}{2}\right \rceil h.$$

\rem \label{extending from half-period} If we know that a word has period $p$ then any factor of length $p$ will determine the rest of the word.  In the case of a palindromic periodicity a factor containing two essential centres will determine the rest of the word.   The maximum alphabet size of a word with period $p$ is $p$.  It's not hard to see that the maximum alphabet size for a palindromic periodicity with half-period $h$ is $h$, $h+1/2$ or $h+1$ if the essential centres are centres of palindromes which are respectively all even, alternately odd and even or all odd.
\section{Creating a Palindromic Periodicity}
Palindromic periodicities arise naturally through the interaction of palindromes and periodicities.  In this section we describe some ways in which this can happen.

\begin{thm}If $w$ is a palindrome which is periodic with period $p$ and $|w| \ge 2p$ then $w$ is a palindromic periodicity with  period $p$.
\begin{proof}
  Say that the palindrome's centre is $c$.
  If the palindrome is even then $w$ is a factor of $v^\omega$ where $v$ is the length $p$ factor that follows $c$.  By periodicity $v$ equals the length $p$ factor immediately preceding $c$ and by palindromicity it equals the reverse of this factor.  Therefore $v$ is a palindrome and $w$ is a palindromic periodicity.

  If both the palindrome and $p$ are odd then $w$ is a factor of $v^\omega$ where $v$ is the length $p$ palindrome centred at $c$.

  If the palindrome is odd and $p$ is even then $w$ is a factor of $(st)^\omega$ where $s$ is the length $p-1$ palindrome with centre $c$ and $t$ is a single letter.
\end{proof}

\end{thm}
\begin{thm} Let $u$ be a finite word, $v$ be its reverse and $|u|=n$. Then if $w$ is a prefix of $u'u^\omega$ and also of $v'v^\omega$,
 where $u'$ is a suffix of u and $v'$ is a suffix of $v$, then $w$ is a palindromic periodicity with period $n$.
\begin{proof}
Let $|v'|=k$ and without loss of generality suppose $u'=\epsilon$. Then $w[k+1..n]$ is a suffix of $u$ and a prefix of $v$. But a length $n-k$ prefix of $v$ is the reverse of a length $n-k$ suffix of $u$ so $w[k+1..n]=R(w[k+1..n])$ and $w[k+1..n]$ is a palindrome with center $(k+1+n)/2$.  Similarly $w[n+1..k+n]$ is a palindrome with center $(2n+k+1)/2$ so that $w[k+1..k+n]$ is the concatenation of two palindromes.  Since $w$ has period $n$ it is a palindromic periodicity with period $n$.
\end{proof}
\end{thm}


\begin{thm} \label{intersecting palindromes} Let $P_1$ and $P_2$ be palindromes with centres $c_1$ and $c_2$ respectively, and radii $r_1$ and $r_2$ respectively in a word $w$ and $c_2>c_1$. Then if  $P_1$ and $P_2$
contain each other's centres and are such that neither is a proper
factor of the other then their union is a palindromic periodicity with period $2(c_2-c_1)$.
\begin{proof} Say the palindromes are in a word $w$ and assume that
\begin{equation}\label{r_2<<r_1}
  r_2 \le r_1.
\end{equation} This involves no loss of generality since if it didn't hold we could take the reverse of $w$. Since the palindromes contain each other's centres we have $c_2-r_2 \le c_1$ and $c_2 \le c_1+r_1$ and
since neither palindrome is a proper factor of the other we have $c_1-r_1 \le c_2-r_2$ and $ c_1+r_1 \le c_2+r_2$. Combining all this gives
\begin{equation} \label{e2}
c_1-r_1 \le c_2-r_2 \le c_1 < c_2 \le c_1+r_1 \le c_2+r_2.
\end{equation} so that $P_1 \cup P_2=w[c_1-r_1\dots c_2+r_2]$. Suppose
\begin{equation}\label{irange}
  i \in [c_1-r_1,c_2+r_2-2(c_2-c_1)]=[c_1-r_1,2c_1+r_2-c_2].
\end{equation}  From
(\ref{e2}) we have $$2c_1+r_2-c_2 \le c_1+r_1$$ so that $w[i]$ is in $P_1$
and by (\ref{e1}) we have
$$w[i]=w[2c_1-i].$$ Now
\begin{align*}2c_1-i & \in
[2c_1-(2c_1+r_2-c_2),2c_1-(c_1-r_1)]\\&=[c_2-r_2, c_1+r_1]\end{align*}
so, by (\ref{e2}), $w[2c_1-i]$ is in $P_2.$ By (\ref{e1})
again we have
$$w[2c_1-i]=w[i+2(c_2-c_1)]$$ for any $i$ satisfying (\ref{irange}) and the maximum value of $i$ maps onto $c_2+r_2$ which is the upper bound of $P_1 \cup P_2$, thus $P_1 \cup P_2$ has period $2(c_2-c_1)$. We now show that $w$ is a palindromic periodicity.\\

  If $c_1$ is an integer then $s=w[c_1]$ and $t=w[c_1+1..2c_2-c_1-1]$ are palindromes with centres $c_1$ and $c_2$ respectively. If $c_1$ is not an integer then $s=w[c_1-1/2..c_1+1/2]$ and
  \begin{eqnarray*}
 t&=&w[c_1+3/2..c_1+3/2+2(c_2-c_1)-3]\\
  &=& w[c_1+3/2..2c_2-c_1-3/2]
 \end{eqnarray*} are palindromes.  In both cases $|s|+|t|=2(c_2-c_1)$ and $s$ is nested in $P_1$ and $t$ is nested in $P_2$. Using (\ref{e2}) it is easily shown that $|P_1 \cup P_2| \ge |s| + |t|$. As noted above, the union of  $P_1$ and $P_2$  has period $2(c_2-c_1)$ and so is a factor of $(st)^{\omega}$, and so is a palindromic periodicity.
\end{proof}
\end{thm}

\begin{rem} In the special case where $c_1=n-1/2$ and $c_2=n+1/2$ for some integer $n$ the period is 1 as well as 2. If $c_1=n$ and $c_2=n+1$ the period is 2 and need not be 1. \end{rem}

 If the condition  of the theorem that neither palindrome  is a proper factor of the other doesn't hold we still have the following result.

\begin{cor} 
Let $w[1..n]$ and $w[k+1,,k+l]$ be palindromes with $1<k+1<k+l<n$, with centers $c_1=(1+n)/2$ and $c_2=k+(l+1)/2$ respectively, that contain each others centres, so that
$ k+1 \le (1+n)/2 \le k+l$. \\
\noindent (a) If $c_1<c_2$ then $w[k+1..n-k]$ is a palindromic periodicity with period $2(c_2-c_1)=n-l-2k$.\\
\noindent (b) If $c_2<c_1$ then $w[n+1-k-l..k+l]$ is a palindromic periodicity with period $2(c_1-c_2)=2k+l-n$.
\begin{proof}
We prove case (a), the other follows by symmetry.  Note that $w[k+1..n-k]$ is a palindrome nested in $w[1..n]$.  The palindrome $w[k+1..k+l$ is a prefix of this so we can apply the theorem and the result follows.
\end{proof}
\end{cor}

We also have the following result when the palindromes don't contain each other's centres.

\begin{thm} \label{thm 2.6}
  If two palindromes intersect, neither is a proper factor of the other and at least one of the two does not contain the centre of the other then their union is a palindromic periodicity whose half-period equals the distance between the palindromes' centres.
  \begin{proof}
    Let the palindromes be $w[a..b]$ and $w[c..d]$ with $\frac{a+b}{2}<\frac{c+d}{2}$.  The conditions of the theorem mean that $a \le c \le b \le d$. Consider the factor $w[c..b]$.  This is a suffix of the first palindrome and so $w[a..d]$ has a prefix equal to $R(w[c..b])$. The factor $w[c..b]$ is also a prefix of the second palindrome so that $w[a..d]$ also has a suffix equal to $R(w(c..b))$.  Thus $w[a..d]$ has a border of length $|w[b-c]|=b-c+1.$  It therefore has period
    $$|w[a..d]|-(b-c+1)=d-a-b+c$$ which is twice the difference between their centres.
      Suppose now that $w[a..b]$ does not contain the centre of $w[c..d]$, so that $b<(c+d)/2$. Then $w$ is a factor of $(ps)^\omega$ where $p$ is the palindrome $w[a..b]$ and $s$ is the palindrome $w[b+1..c+d-b-1]$ which is nested in $w[c..d]$.
  \end{proof}
\end{thm}

\section{Towards a periodicity lemma for palindromic periodicities}

 We say that a word is a \emph{double palindromic periodicity} if it is a palindromic periodicity in two ways each with its own offset and half-period. In this section we obtain something like Fine and Wilf's lemma but applied to words which are double palindromic periodicities.  That is, we show that if a double palindromic periodicity with half-periods $h_1$ and $h_2$ is sufficiently long then it is a palindromic periodicity whose half-period is possibly less than both $h_1$ and $h_2$.

 To do this we first consider the case of a single palindrome embedded in a single  palindromic periodicity and show that if the embedded  palindrome is sufficiently long the half-period of the palindromic periodicity is reduced.  Here ``sufficiently long" depends on the parameters of the palindromic periodicity and on the position of the centre of the embedded palindrome.  Some preliminary results are needed  to prove this.  They are Lemma  \ref{3.2}, Theorem \ref{3.3} and Corollary \ref{3.4}. The main theorem about an embedded palindrome is Theorem \ref{g-word has period g}.
 A palindromic periodicity contains palindromic prefixes and suffixes.  If one of these is sufficiently long compared to $h_1$ or $h_2$, say to $h_1$, we can apply Theorem \ref{g-word has period g} and the whole double palindromic prefix takes on a half-period, $h_3$ say, which is less than $h_1$.  The word now has periods $2h_3$ and $2h_2$.  Finally we apply Fine and Wilf's Lemma to obtain our main theorem.

 We now define something called a $g$-word.  This will be the palindrome we embed in our double palindromic periodicity.

A \emph{g-word} is a non-empty palindrome $w[i+1..i+n]$ with centre $c$ which is embedded in  a palindromic periodicity $w$ with half-period $h$ and offset $r$.  Its length is $n=2h-g$ where $$g=\gcd(2|c-r|,2h).$$ Saying that the palindrome is non-empty means we do not allow $2h$ to equal $g$.  For example,

\begin{large}
\begin{center}
\begin{verbatim}
        abcdexxedcbaabcdexxedcbaabcdexxedcbaabcdexxedcba
\end{verbatim}
\end{center}
\end{large}
 is a g-word with parameters $r=13/2$, $c=49/2$, $h=30$, $n=48$ and $g=12$.

\begin{rem} Notice that the length of a $g$-word is always divisible by $g$. \end{rem}

The next results show that a $g$-word is a palindromic periodicity with period $g$. Clearly the word above does have period $g=12$.
The following lemma tells us a bit about the structure of a $g$-word.

\begin{lem}\label{3.2} A g-word $w[i+1..i+n]$ contains exactly two essential centres of the palindromic periodicity, $i+r$ and $i+r+h$.

\begin{proof} In this proof we set $i=0$ which does not cause any loss of generality.  We know that $w$ contains an essential palindromic centre at $r$ and, since its length is less than $2h$, it can contain at most one other, and this must be at $r+h$. Clearly $w[1..2h]$ contains exactly two such centres. We must show that $r+h$ does not lie in the factor $w[n+1..2h]$, that is, that $r+h\le n+1/2$.

We have \begin{eqnarray*}
g &=& gcd(|2c-2r|,2h) \\
 n &=& 2h-g \\
  c &=& (n+1)/2.
  \end{eqnarray*} In terms of $c$, $g$ and $r$,
  \begin{eqnarray*}
    n &=& 2c-1 \\
    h &=& c+(g-1)/2\\
    g &=& \gcd(|2c-2r|,2c+g-1)\\
      &=& \gcd(|2c-2r|,2c-1).
  \end{eqnarray*}
 We consider two cases.

 If $r \le c$ then $g=\gcd(2c-2r,2c-1)$.  Suppose, for the sake of contradiction, that  $r+h>n+1/2$.  Then
 \begin{equation*}
   r+c+(g-1)/2 > 2c-1/2 \;
    \Rightarrow  \; g>2c-2r.\\
  \end{equation*}
  But by definition $g$ divides $2c-2r$ so we have a contradiction and conclude that if $c \ge r$ then $r+h \le n+1/2$, and so $r+h$ does not lie in $w[n+1..2h]$.

 On the other hand, if $r>c$ then $g=\gcd(2r-2c,2c-1)$. Suppose, for the sake of contradiction, that $r+h \le 2h$. Then
 \begin{eqnarray*}
 r \le h &\Rightarrow& r \le c+(g-1)/2 \\
   &\Rightarrow& 2r-2c \le g-1\\
   & \Rightarrow & g>2r-2c.
 \end{eqnarray*} But this contradicts our definition which requires that $g$ divides $2r-2c$ so we conclude that if $r> c$ then $r+h>2h$ so $r+h$ does not lie in $w[n+1..2h]$ and we are done.
\end{proof}
\end{lem}
We have five parameters to use in describing a $g$-word: $r$, $c$, $h$, $n$ and $g$.  These are not independent.  We can determine all five if we know $c-r$ and $h$. However we'll find it convenient to use all of them. It's possible that either $r$ or $r+h$ coincides with $c$. In this case either $g =\gcd(0,2h)=2h$ or $r+h=c$ andl$g=\gcd(2h,2h)=2h$. In either case $n=2h-g=0$ and $w$ is empty. We'll assume henceforth that neither $r$ nor $r+h$ coincides with $c$.  In this case we have $r<c<r+h$.

\begin{thm}\label{3.3}
  (a) A $g$-word is a power of a length $g$ palindrome, and\\
  (b) this periodicity extends to the whole of the embedding palindromic periodicity.
 
  \begin{proof}
    (a) Let $w$ be a $g$-word as in Lemma \ref{3.2} and without loss of generality set $i=0$ so that $w=w[1..2h-g]$ with essential centers at $r$ and $r+h$. Thus $w[1..2r-1]$ and $w[2(h+r)-(2h-g)..2h-g]=w[2r+g..2h-g]$ are, respectively, a palindromic prefix and a palindromic suffix of $w$. Since a $g$-word is a palindrome $w$ also has a palindromic centre at $c$.
    
    The palindrome centred at $r$ is contained in $w$ so we may apply Theorem \ref{intersecting palindromes} or Theorem \ref{thm 2.6} (depending on whether or not $c$ is inside $w[1..2r-1]$) and conclude that $w$ has period $2(c-r)$. Similarly, using the palindrome centred at $r+h$, we find that $w$ also has period $2(h+r-c)$.  We now apply the Fine Wilf lemma. The greatest common divisor of these periods is
    \begin{equation*}
      (2c-2r,2h+2r-2c = (2(c-r),2h)=g \\   
     \end{equation*} and their sum is $2h$. Since $|w|=2h-g$ we may apply the lemma and conclude that $w$ has period $g$.  Since the length $g$ prefix equals the reverse of the length $g$ suffix and $|w|$ is divisible by $g$ this prefix is a palindrome $w$ is a power of it.\\
     (b) The embedding palindromic periodicity has period $2h$, the $g$-word has period $g$ which divides $2h$ so part (b) follows from Lemma \ref{L4}.
\end{proof}

\end{thm}

The following example shows that the theorem is sharp in the sense that it would not hold if a $g$-word was defined to have some length less than $2h-g$.\\

\noindent \textbf{Example}. The first word below is a palindromic periodicity with offset 25/2, half-period  30 and length 60. In the second word we have inserted a $g$-word at the centre of this palindromic periodicity whose parameters are $r=25/2$, $c=49/2$, $g=12$ and length is $48$. This word has period 12 in agreement with Theorem \ref{g-word has period g} and the periodicity extends to the whole of the word. In the last word the central palindrome has length 46 which is two letters too short to be a $g$-word. In this case the central palindrome has period 24. This periodicity does not extend to the whole word.  This shows that in some cases Theorem \ref{3.3} is sharp.\\

\begin{large}
\begin{center}
\begin{verbatim}
 abcdefghijkllkjihgfedcbamnopqrstuvwxyzABCDDCBAzyxwvutsrqponm
 abcdeffedcbaabcdeffedcbaabcdeffedcbaabcdeffedcbaabcdeffedcba
 cdefgbbgfedccdefgbbgfedccdefgaagfedccdefgbbgfedccdefgaagfedc
\end{verbatim}
\end{center}
\end{large}

We now present our periodicity theorem for palindromic periodicities.

 \begin{thm}\label{periodicity lemma} If $w$ is a double palindromic periodicity with parameters $(r_1,h_1)$ and $(r_2,h_2)$ and with
\begin{equation}\label{length of dpp}
|w|\ge 2h_1+2h_2-\gcd(2(r_2-r_1),2h_1,2h_2)
\end{equation} then $w$ has period $\gcd(2(r_2-r_1),2h_1,2h_2)$.
\begin{proof}
We write $g$ for $\gcd(2(r_2-r_1),2h_1,2h_2)$ and let $c=(|w|+1)/2$ be the centre of $w$. By (\ref{length of dpp}) we have

  $$2c-1 \ge 2h_1+2h_2-g $$
so that $$r_1 + (\frac{c-r_1}{h_1}-1)h_1\ge h_2-\frac{g-1}{2}.$$
 But $$\frac{c-r_1}{h_1}-1 < \lfloor \frac{c-r_1}{h_1}\rfloor$$ so $$r_1+\lfloor \frac{c-r_1}{h_1}\rfloor h_1 > h_2-\frac{g-1}{2}.$$
The left hand side here is the position of largest essential centre of the first palindromic periodicity which is less than or equal to $c$.  Writing $i$ for $\lfloor \frac{c_1-r_1}{h_1}\rfloor$ we get
$$2(r_1+ih_1)-1 > 2h_2-g.$$
The left hand side  is the length of the palindromic prefix of $w$ centred at $r_1+ih_1$. Note that $g$ divides $\gcd(2h_2,2(r_1+ih_1-r_2)$ so that
$$2(r_1+ih_1)-1 \ge 2h_2-\gcd(2h_2,2(r_1+ih_1-r_2))$$
 and we see that this palindromic prefix contains a $g$-word. Then Theorem \ref{g-word has period g} applies and $w[1..2(r_1+ih_1)-1]$ has period $\gcd(2(r_2-(r_1+i)h_1),2h_2)$ and this periodicity extends to the whole of $w$. But $w$  also has period $2h_1$. In order to apply Fine and Wilf's  Lemma \ref{FW} we must show that the length of $w$ is at least \begin{eqnarray*}
&&2h_1+\gcd(2(r_2-(r_1+i)h_1),2h_2) -\gcd(2(r_2-(r_1+i)h_1),2h_2 ,2h_1)\\
                                            &=& 2h_1+\gcd(2(r_2-(r_1+i)h_1),2h_2) -g
                                         \end{eqnarray*}   which is equivalent to showing that $2h_2 \ge \gcd(2(r_2-(r_1+i)h_1),2h_2)$ which is clearly true. So $w$ has period $$\gcd(2(r_2-(r_1+ih_1)),2h_1,2h_2)=\gcd(2(r_2-r_1),2h_1,2h_2)$$ as required.
\end{proof}
\end{thm}
Since $gcd(2(r_2-r_1),2h_1,2h_2)$ may be less $\gcd(2h_1,2h_2)$ our theorem may need a longer word to apply than does the Fine and Wilf Lemma.  On the other hand our result may imply a shorter period than does Fine and Wilf.

Note that if $2(r_2-r_1)$ divides $h_1$ and $h_2$ then the theorem turns into the Fine and Wilf Lemma.
Unlike Fine and Wilf's result, this theorem is  usually not sharp.  That is, for many combinations of $h_1$, $h_2$, $r_1$ and $r_2$ words shorter than $2h_1+2h_2-\gcd(2(r_2-r_1),2h_1,2h_2)$ will have period $\gcd(2(r_2-r_1),2h_1,2h_2)$.  This is illustrated in Table 1 which shows the least periods of double palindromic periodicities with parameters $h_1=4$, $h_2=6$ and all combinations of $r_1$ and $r_2$ for which $r_1$ and $r_2$ have the same  parity. The last stipulation means that $2(r_2-r_1)$ is divisible by 4 and $\gcd(2(r_2-r_1),2h_1,2h_2)=4$. The table shows the least period of words of decreasing length starting with length $2h_1+2h_2-\gcd(2(r_2-r_1),2h_1,2h_2)=16$. A worthwhile periodicity lemma would predict all periods in this table, not just those in the third column. In Table 2 $r_1$ and $r_2$  have opposite parity so that $\gcd(2(r_2-r_1),2h_1,2h_2)=2$ and lengths are decreasing from 18.

    \begin{figure}
  \begin{tabular}{cclccccccccccc}
 $r_1$ & $r_2$ &\text{\emph{lengths}}&16&15&14&13&12&11&10&9&8&7&6\\ \hline
 0&0&&4&4&4&4&4&4&4&4&4&4&6\\
\noalign{\medskip}0&2&&4&8&8&8&8&8&8&8&8&4&4\\
\noalign{\medskip}0&4&&4&4&4&4&4&8&8&8&8&6&6\\
\noalign{\medskip}1&1&&4&4&4&4&4&4&4&4&8&7&6\\
\noalign{\medskip}1&3&&4&4&4&4&4&4&4&4&4&4&4\\
\noalign{\medskip}1&5&&4&4&4&4&8&8&8&8&8&7&6\\
\noalign{\medskip}2&0&&4&4&4&8&8&8&8&8&8&7&6\\
\noalign{\medskip}2&2&&4&4&4&4&4&4&4&8&8&7&6\\
\noalign{\medskip}2&4&&4&4&4&4&4&4&4&4&4&4&4\\
\noalign{\medskip}3&1&&4&4&8&8&8&8&8&8&8&7&6\\
\noalign{\medskip}3&3&&4&4&4&4&4&4&8&8&8&7&6\\
\noalign{\medskip}3&5&&4&4&4&4&4&4&4&4&4&4&4
\end {tabular}
\caption {Table 1. Periods of double palindromic periodicities with parameters $h_1=4$, $h_2=6$ and offsets and lengths as shown. In each case $\gcd(2(r_2-r_1),2h_1,2h_2)=4$.  }

\end{figure}
 \begin{figure}
 \begin{tabular}{cclccccccccccc}
 $r_1$ & $r_2$ &\text{\emph{lengths}}&18&17&16&15&14&13&12&11&10&9&8\\ \hline
 0&1&&2&2&2&2&2&2&2&2&2&6&6\\
 \noalign{\medskip}0&3&&2&2&2&2&2&2&2&2&2&8&8\\
 \noalign{\medskip}0&5&&2&2&2&2&2&2&2&2&2&2&2\\
\noalign{\medskip}1&0&&2&2&2&2&2&2&2&2&2&2&8\\
\noalign{\medskip}1&2&&2&2&2&2&2&2&2&2&6&6&6\\
\noalign{\medskip}1&4&&2&2&2&2&2&2&2&2&2&2&2\\
\noalign{\medskip}2&1&&2&2&2&2&2&2&2&2&2&8&8\\
\noalign{\medskip}2&3&&2&2&2&2&2&2&2&6&6&6&6\\
\noalign{\medskip}2&5&&2&2&2&2&2&2&2&2&2&2&2\\
\noalign{\medskip}3&0&&2&2&2&2&2&2&2&2&2&2&6\\
\noalign{\medskip}3&2&&2&2&2&2&2&2&2&2&2&8&8\\
\noalign{\medskip}3&4&&2&2&2&2&2&2&2&2&2&2&2
\end {tabular}

\caption {Table 2. Periods of double palindromic periodicities with parameters $h_1=4$, $h_2=6$ and offsets and lengths as shown. In each case $\gcd(2(r_2-r_1),2h_1,2h_2)=2$. Here all words of length at least 12 have period 2.}

\end{figure}

\section{Discussion} There are several directions in which research into palindromic periodicities might proceed.  One is to obtain a stronger version of Theorem \ref{periodicity lemma}. Another would be to look for and count occurrences of palindromic periodicities in famous words such as the Thue-Morse word, Sturmian words, particularly the Fibonacci word, and the Oldenburger-Kolakoski word. This may require a clever way of recognising and counting palindromic periodic factors.

My interest in palindromic periodicities began with the paper \cite{JS} in which an almost sharp bound was obtained for the maximum number of distinct palindrome in circular words.  A major ingredient in that paper was a weaker version of Theorem \ref{intersecting palindromes}.  This version recognised that the union of a pair of palindromes containing  each other's centres is periodic, but not that it was a palindromic periodicity.  It might be that Theorem \ref{intersecting palindromes}  might lead to a sharp bound. Similarly paper \cite{GSS} by Glen, Simpson and Smyth failed to obtain a sharp bound  on the maximum  number of distinct palindromes in an edge-labelled starlike tree.

Another direction would be the investigate the maximum number of maximal palindromic periodicities in a word.  This would be analogous to the problem of determining the maximum number of maximal periodicities, also know as runs, that can occur in a word of length $n$.  These are periodic factors, with length at least twice the period, which cannot be extended to the left or right without altering their periods. In 2000 Kolpakov and Kucherov \cite{KK} showed that the number was $O(n)$, without giving any information about the size of the implied constant. They conjectured that the number of runs was less than $n$.  Then Rytter \cite{R} showed the number was less than $5n$.  This was followed by a sequence of increasingly long and complicated papers decreasing the bound.  Then along came Bannai, I, Inenaga, Nakashima, Takeda, and  Tsuruta who showed, with a very short and elegant proof, that the Kolpakov-Kucherov conjecture was correct.  Since then the bound has been further decreased, best so far being $183/193$ by \u{S}t\u{e}p\'{a}n Holub \cite{SH}.

\bibliography{}

\end{document}